\documentclass[11pt]{article}
\usepackage{amsfonts}
\usepackage{fancyhdr}
\usepackage[dvips, usenames]{color}

\usepackage{latexsym}
\usepackage{amsmath}
\usepackage{amssymb}
\usepackage{epsfig}
\usepackage{mathrsfs}
\usepackage{cite}

\usepackage{lineno}

\begin{document}

\newtheorem{thm}{Theorem}[section]
\newtheorem{cor}{Corollary}[section]
\newtheorem{lem}{Lemma}[section]
\newtheorem{definition}{Definition}[section]
\newtheorem{con}{Conjecture}[section]

\begin{center}
{\Large\bf The super-connectivity of Kneser graph KG(n,3)}
\\[20pt]
{Yulan\ Chen $^a$ 
\quad Yuqing\ Lin $^b$ \footnote{Corresponding Author.} \quad Weigen\
Yan $^a$
\footnote{Partially supported by NSFC Grant (12071180; 11571139).\\
{\it Email address:} y\_l\_Chen@163.com (Y. Chen),
 Yuqing.lin@newcastle.edu.au, weigenyan@263.net (W. Yan)}}
\\[5pt]
\footnotesize {$^a$ School of Sciences, Jimei University,
Xiamen 361021, China \\
$^b$ School of Electrical Engineering, the University of Newcastle, Australia}

\end{center}
\begin{abstract}
A vertex cut $S$ of a connected graph $G$ is a subset of vertices of $G$ whose deletion makes $G$ disconnected. A super vertex cut $S$ of a connected graph $G$ is a subset of vertices of $G$ whose deletion makes $G$ disconnected and there is no isolated vertex in each component of $G-S$. The super-connectivity of graph $G$ is the size of the minimum super vertex cut of $G$. Let $KG(n,k)$ be the Kneser graph whose vertices set are the $k$-subsets of $\{1,\cdots,n\}$, where $k$ is the number of labels of each vertex in $G$. We aim to show that the conjecture from Boruzanli and Gauci \cite{EG19} on the super-connectivity of Kneser graph $KG(n,k)$ is true when $k=3$. 
\end{abstract}

\section{Introduction} \

Let $[n]=\{1,\cdots,n\}$ be $n$ labels, the Kneser graph $G=KG(n,k)$ is the graph whose vertices are the $k$-subsets of $[n]$, two vertices are adjacent if these two $k$-subsets are disjoint, i.e. two vertices do not share labels. Let $V(G)$ be the set of vertices of $G$, it is clear $V(KG(n,k)) = {[n] \choose k}$ and the $KG(n,k)$ is regular with degree ${n-k \choose k}$. A vertex cut $S$ of a connected graph $G$ is a subset of vertices of $G$ whose deletion disconnect $G$. The connectivity $\kappa$ of $G$ is the size of the minimum vertex cut of $G$. If the deletion of any vertex cut of size $\kappa$ in $G$ will isolate a vertex, then $G$ is super-connected. A vertex cut which isolate a single vertex is called an trivial vertex cut of $G$. When $G$ is super-connected, it makes sense to determine the size of minimum nontrivial vertex cut of $G$, that is, the super-connectivity $\kappa_{1}$ of $G$. And the smallest nontrivial vertex cut is called a super-vertex cut of $G$. A complete graph $K_n$ is a simple graph with $n$ vertices and edge between every pair of vertices of $K_n$. 

The concept of Kneser graph was proposed by Kneser in 1955 \cite{Kneser55}. Structural properties of Kneser graph has been studied extensively, for example, the hamiltoniancity, Chromatic number and the matchings in Kneser graph. Chen and Lih proved that the Kneser graph is symmetric, vertex-transitive and edge-transitive \cite{CL87}. Using this property, Boruzanli and Gauci showed that the connectivity of Kneser graph $KG(n,k)$ is ${n-k \choose k}$ \cite{EG19}. Harary has proposed the concept of super-connectivity in 1983 \cite{Harary83}. Subsequently, Balbuena, Marcote and Garc\'{\i}a-V\'{a}zquez defined a similar concept, i.e. restricted connectivity of graphs \cite{BMG05}. In this paper, we will investigate the super-connectivity of Kneser graph. 

It is known that if $n < 2k$, then $KG(n,k)$ contains no edges, and if $n = 2k$, then $KG(n,k)$ is a set of independent edges. The Kneser graph $KG(n, 1)$ is the complete graph on $n$ vertices. Boruzanli and Gauci has made a conjecture which states that \cite{EG19},
\begin{con}
 Let $n\geq2k+1$, then the super-connectivity $\kappa_{1}$ of $KG(n,k)$ is
\[
   \kappa_{1}=\left\{\begin{array}{lcl}
                       2\left({n-k \choose k}-1 \right) & if & 2k+1\leq n<3k, \\
                       2\left({n-k \choose k}-1 \right)-{n-2k\choose k} & if & n \geq 3k.
                     \end{array}
   \right.
\]
\end{con}

Boruzanli and Gauci has proved that this conjecture holds when $k=2$\cite{EG19}, in this work, we are considering the case when $k=3$.

\section{Super-Connectivity of $KG(n,3)$}

In this section, we determine the super-connectivity of $KG(n,3)$ when $n \geq 7$ and confirmed that the Conjecture 1.1 is true for $k=3$. 

\begin{thm}
 The super-connectivity of Kneser graph $KG(n,3)$ is
 \[
   \kappa_{1}=\left\{\begin{array}{lcl}
                       2\left({n-3 \choose 3}-1 \right) & if & 7\leq n \leq 8, \\
                       2\left({n-3 \choose 3}-1 \right)-{n-6\choose 3} & if & n \geq 9.
                     \end{array}
   \right.
 \]
\end{thm}

\textit{Proof}. 

Let $S\subseteq V(G)$ be a super-vertex cut of $G$. Suppose $n \geq 9$ and  $|S| <  2\left({n-3 \choose 3}-1 \right)-{n-6\choose 3}$, then we have
\begin{align*}
   |G-S| & > {n\choose 3}-2\left[{n-3\choose 3}-1\right]+{n-6 \choose 3} = \frac{54n-204}{6}=9n-34
\end{align*}

This means that if the $\kappa_{1}$ is less than the bound stated in the conjecture, then there will be more than $9n-34$ vertices in $G-S$. In the following, we will show that $G-S$ has to be connected if it has more than $9n-34$ vertices.

Since $S$ is a super-vertex cut, then $G-S$ has at least two components and each component has at least $2$ vertices. If $G-S$ has a component containing only two vertices, then it is straightforward that $|S| = \kappa_{1}$ since the $S$ has to contain all the neighbours of these two vertices and also it is easy to see that there is no singular vertex in $G-S$.

Now we assume that each component of $G-S$ has at least $3$ vertices. We also assumed that $G-S$ has two components $C_{1},C_{2}$, and $C_2=G-S-C_1$. Note, in here, $C_2$ might not be connected. If $C_{2}$ is not connected, then $C_{2}$ is the union of some connected components with have at least $3$ vertices in each. Since $C_1$ has at least three vertices, let them be $v_{1},v_{2},v_{3}$. These three vertices form possibly two different graphs, either a complete graph $K_3$ or a path $P_3$ of length $2$, if these three vertices form a path, then there are two possibilities, either the two non-adjacent vertices share only one common label, which we refer to as $Type\ 1 \ path$ or the two non-adjacent vertices share two common labels, which we refer to as $Type \ 2 \ path$.

We are making the following three claims.

Claim 1: If there is a $K_3$ in $C_1$, then there are at most $27$ vertices in $C_2$.

Let the three vertices in $C_1$ be $v_{1}=\{1,2,3\}$, $v_{2}=\{4,5,6\}$, $v_{3}=\{7,8,9\}$. Since $C_1$ and $C_2$ are disconnected, then every vertex in $C_{2}$ has at least one label in common with every vertex in $C_{1}$, i.e. any vertex of $C_2$ has to have a label from $\{1,2,3\}$, a label from $\{4,5,6\}$ and a label from $\{7,8,9\}$, then the number of vertices in $C_{2}$ is at most $3^{3}=27$.

Claim 2: If there is a Type 1 path in $C_1$, then there are at most $3n+3$ vertices in $C_2$.

Let the three vertices in $C_1$ be $v_{1}=\{1,2,3\}$, $v_{2}=\{4,5,6\}$, $v_{3}=\{1,7,8\}$, the common label of the two end vertices is $1$, then similar to the proof of Claim 1, we have maximum $3(n-2)$ vertices in $C_{2}$ contain label $1$, since the vertices of $C_2$ in this case has to use a label in $\{4,5,6\}$, and in this calculation we have double counted $3$ vertices $\{1,4,5\}$, $\{1,4,6\}$, $\{1,5,6\}$, therefore, there are at most $3(n-2)-3$ vertices containing label $1$ in $C_2$. And there are at most $2\cdot3\cdot2$ vertices in $C_{2}$ do not contain label $1$. Hence the number of vertices in $C_{2}$ is at most $3(n-2)-3+12=3n+3$.

Claim 3: If there is a Type 2 path $P_3$ in $C_1$, then there are at most $6n-18$ vertices in $C_2$.

Let the three vertices in $C_1$ be $v_{1}=\{1,2,3\}$, $v_{2}=\{4,5,6\}$, $v_{3}=\{1,2,7\}$, the set of common labels of the end vertices are $\{1,2\}$. Similar to the previous argument, we have maximum $3(n-3)$ vertices in $C_{2}$ contain label $1$, but not label $2$. Similarly, have maximum $3(n-3)$ vertices in $C_{2}$ contain label $2$, but not label $1$. And there are at most $3$ vertices in $C_{2}$ contain both labels $\{1,2\}$, and there are at most $3$ vertices in $C_{2}$ contain neither label $1$ nor label $2$. Since we have double counted the $6$ vertices $\{1,4,5\}$, $\{1,4,6\}$, $\{1,5,6\}$, $\{2,4,5\}$, $\{2,4,6\}$, $\{2,5,6\}$, then the number of vertices in $C_{2}$ is at most $2\cdot3(n-3)+6-6=6n-18$.

Next we will show that $C_1 \cup C_2 \leq 9n-34$, which implies that $G-S$ has to be connected if it contains more than $9n-34$ vertices. We consider the following cases.

Case 1: There are $K_3$s in both $C_1$ and $C_2$.

\begin{figure}[htbp]
  \centering
  \includegraphics[height=6cm]{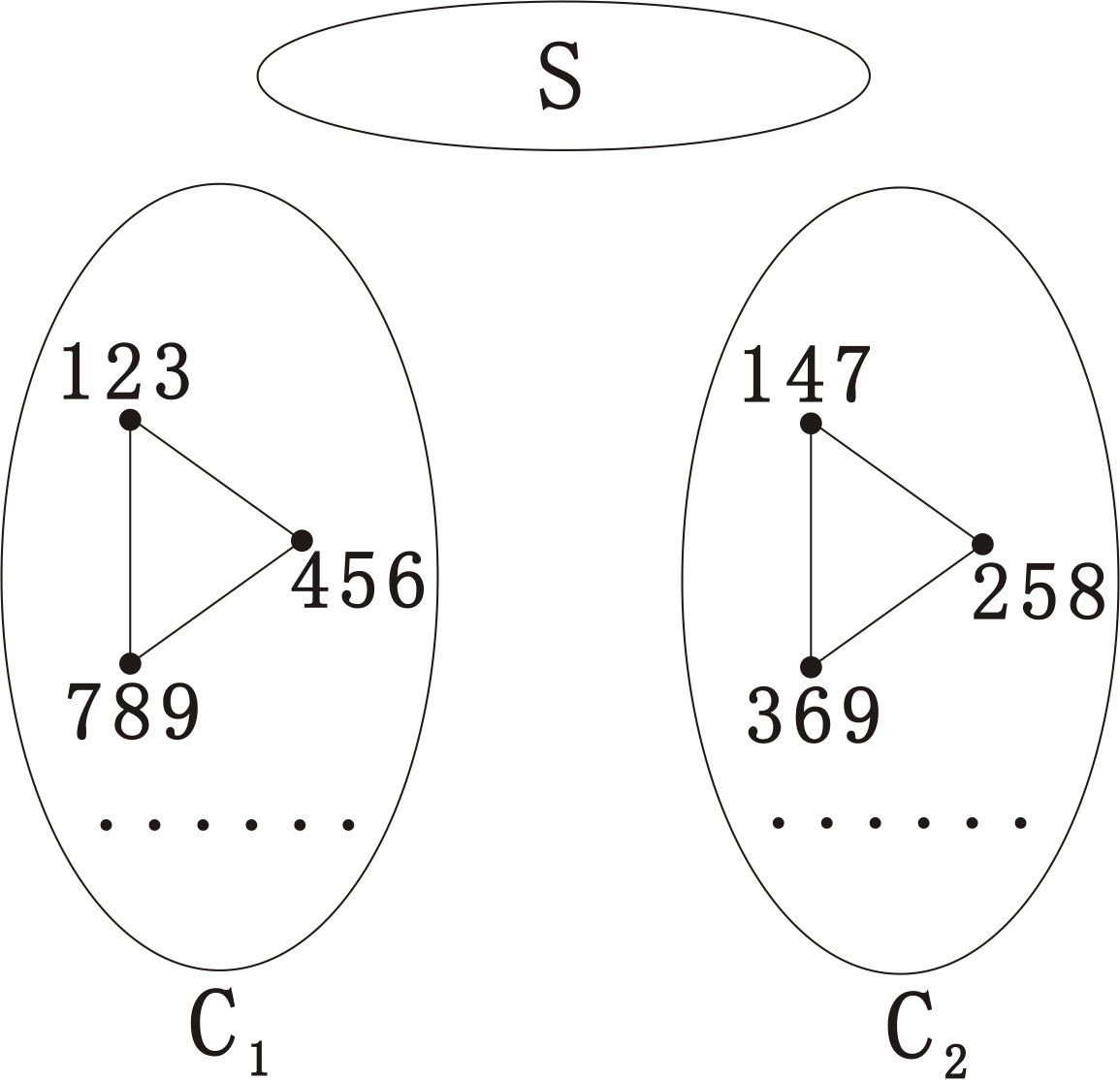}\\
  \caption{The case in $C_1$ and $C_2$}\label{case1}
\end{figure}

If the three vertices in $C_{1}$ form a complete graph $K_3$, let them be $v_{1}=\{1,2,3\}$, $v_{2}=\{4,5,6\}$, $v_{3}=\{7,8,9\}$, then by Claim 1, we have the number of vertices in $C_{2}$ is at most $27$. If $C_2$ also contains a $K_3$, for example, $\{1,4,7\}$, $\{2,5,8\}$, $\{3,6,9\}$, then the $C_1 \cup C_2$ has at most $54$ vertices. In this $54$ vertices, we have double counted the $6$ vertices $\{1,5,9\}$, $\{1,6,8\}$, $\{2,4,9\}$, $\{2,6,7\}$, $\{3,4,8\}$, $\{3,5,7\}$. And the vertex $\{2,3,7\}$ can only be in $C_1$ or $S$, $\{1,4,8\}$ can only be in $C_2$ or $S$, however they are connected, so one of them must be in $S$, the same for the pairs $\{5,6,7\}$ and $\{1,4,9\}$, $\{2,3,4\}$ and $\{1,5,7\}$, thus $C_1 \cup C_2$ has at most $45$ vertices. 
When $n \geq 9$, $9n-34$ is larger than $45$, thus $G-S$ be connected, i.e. $C_1$ and $C_2$ must be connected in this case, a contradiction. So we know that $C_1$ and $C_2$ can not contain $K_3$ at the same time. See Figure \ref{case1} as an example.

Case 2: There is a $K_3$ in $C_1$ or $C_2$, but not in both.

\begin{figure}[htbp]
  \centering
  \includegraphics[height=6cm]{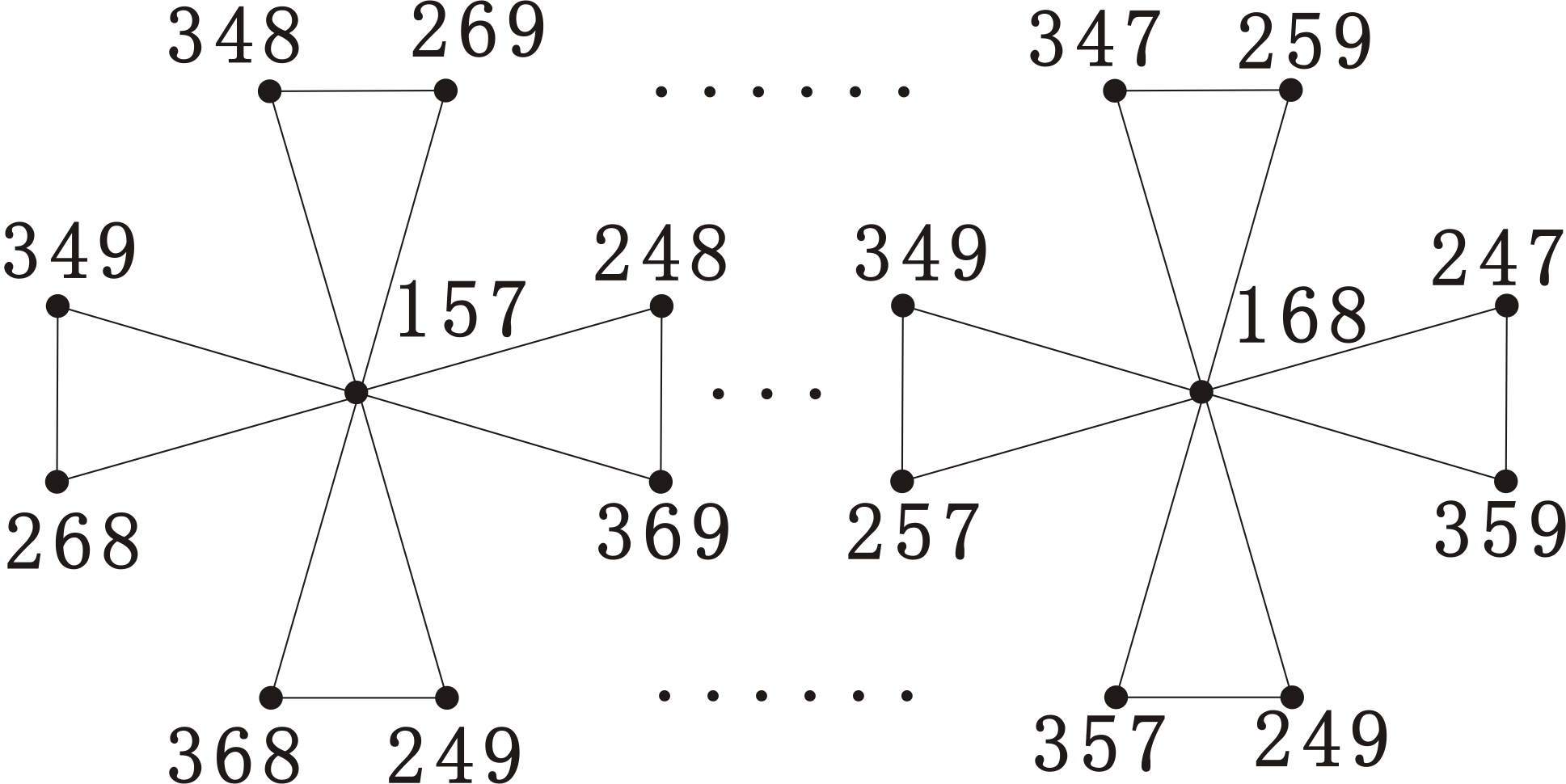}\\
  \caption{The $K_{3}$s in $C_2$}\label{K3s}
\end{figure}

Suppose $C_1$ contains a $K_3$ and $C_2$ do not have a $K_3$. Let the three vertices in $C_1$ be $v_{1}=\{1,2,3\}$, $v_{2}=\{4,5,6\}$, $v_{3}=\{7,8,9\}$. From the Claim 1, we know that there are at most $27$ vertices in $C_2$. If all $27$ vertices are presented in $C_2$, it is easy to verify that there are $36$ $K_3$s in $C_2$, and no two $K_3$s share an edge, however, 4 $K_3$s will share a vertex, for example, $\{1,5,7\}$, $\{2,4,8\}$, $\{3,6,9\}$, and $\{1,5,7\}$, $\{2,4,9\}$, $\{3,6,8\}$, and $\{1,5,7\}$, $\{2,6,8\}$, $\{3,4,9\}$,  and $\{1,5,7\}$, $\{2,6,9\}$, $\{3,4,8\}$ (see Figure \ref{K3s}). It is straightforward to see that least 9 vertices have to be excluded from these $27$ vertices, so that there are no $K_3$ in $C_2$, thus, there are at most $27-9=18$ vertices in $C_2$. If there is exactly $18$ vertices in $C_2$, it implies that these 9 vertices that been removed all contains a certain label, for example, label 3. Otherwise, more than 9 vertices have to be excluded so there is no $K_3$s in $C_2$. 

There must be a path $P_3$ in $C_2$, either Type 1 or Type 2, otherwise, there will be an isolated vertex or $K_2$ in $C_2$, which contradicts the assumption that the number of vertices in each component of $C_2$ is at least $3$.

\begin{figure}[htbp]
  \centering
  \includegraphics[height=6cm]{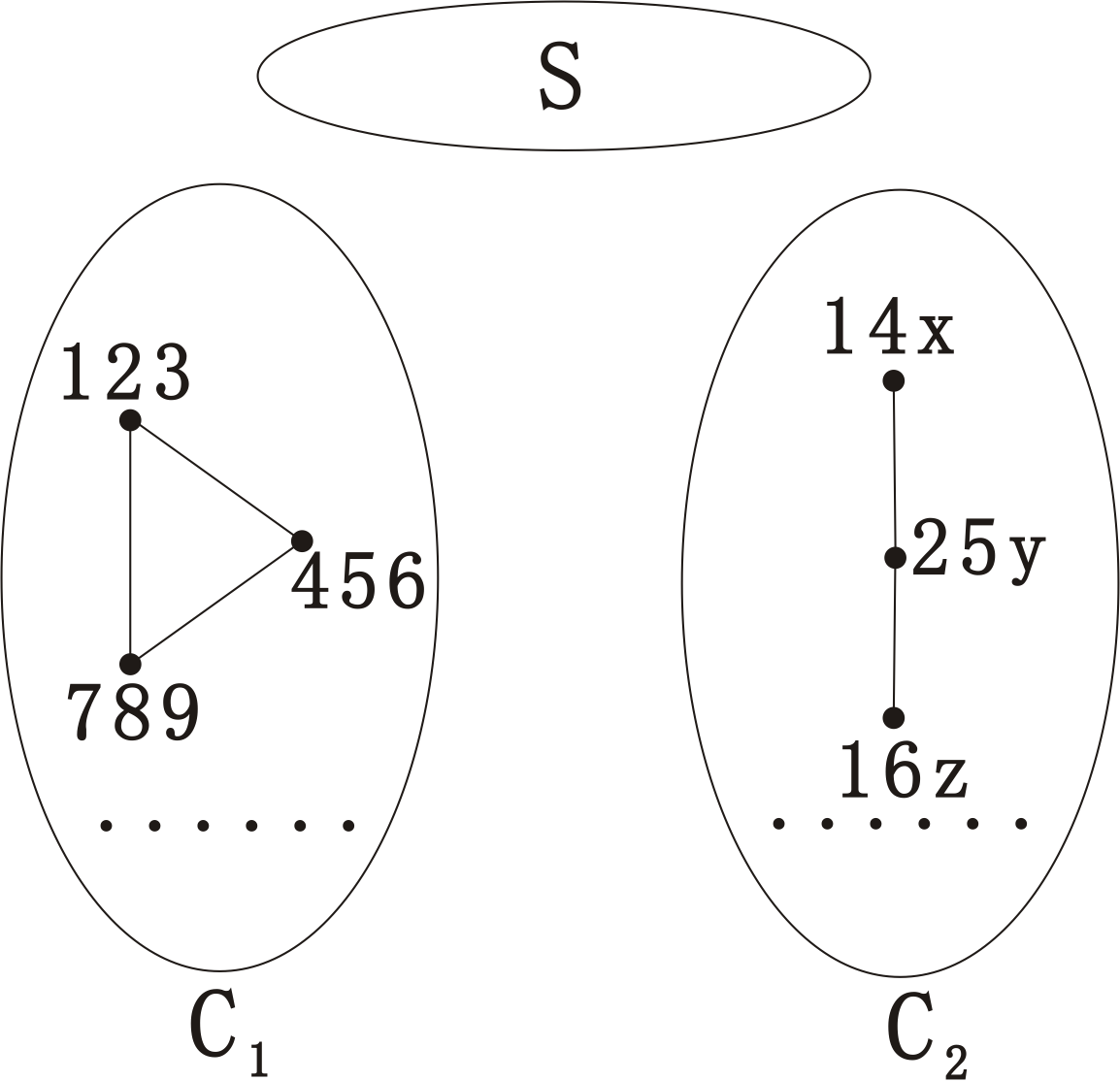}\\
  \caption{The case in $C_1$ and $C_2$}\label{case2.1}
\end{figure}

For the first case, without lose of generality, assume the common label for two end vertices of the path is $1$, and the middle vertex in the path contains label $2$. We could further assume that the three vertices on the path are $\{1,4,x\}$, $\{2,5,y\}$, $\{1,6,z\} \in C_2$, where $x\neq y\neq z$ and $x,y,z \in \{7,8,9\}$. From the proof of Claim 2, we know that there are at most $3n+3$ vertices in $C_1$. However, we have double counted $7$ vertices $\{1,4,y\}$, $\{1,5,7\}$, $\{1,5,8\}$, $\{1,5,9\}$, $\{1,6,y\}$, $\{2,4,z\}$, $\{2,6,x\}$, which should be in either $C_1$ or $C_2$ but not in both. Thus, overall, $C_1 \cup C_2$ has no more than $18+3n+3-7=3n+14$ vertices, which is less than $9n-34$ when $n \geq 9$, then $C_1$ and $C_2$ has to be connected. A contradiction. See Figure \ref{case2.1} for an illustration.

\begin{figure}[htbp]
  \centering
  \includegraphics[height=6cm]{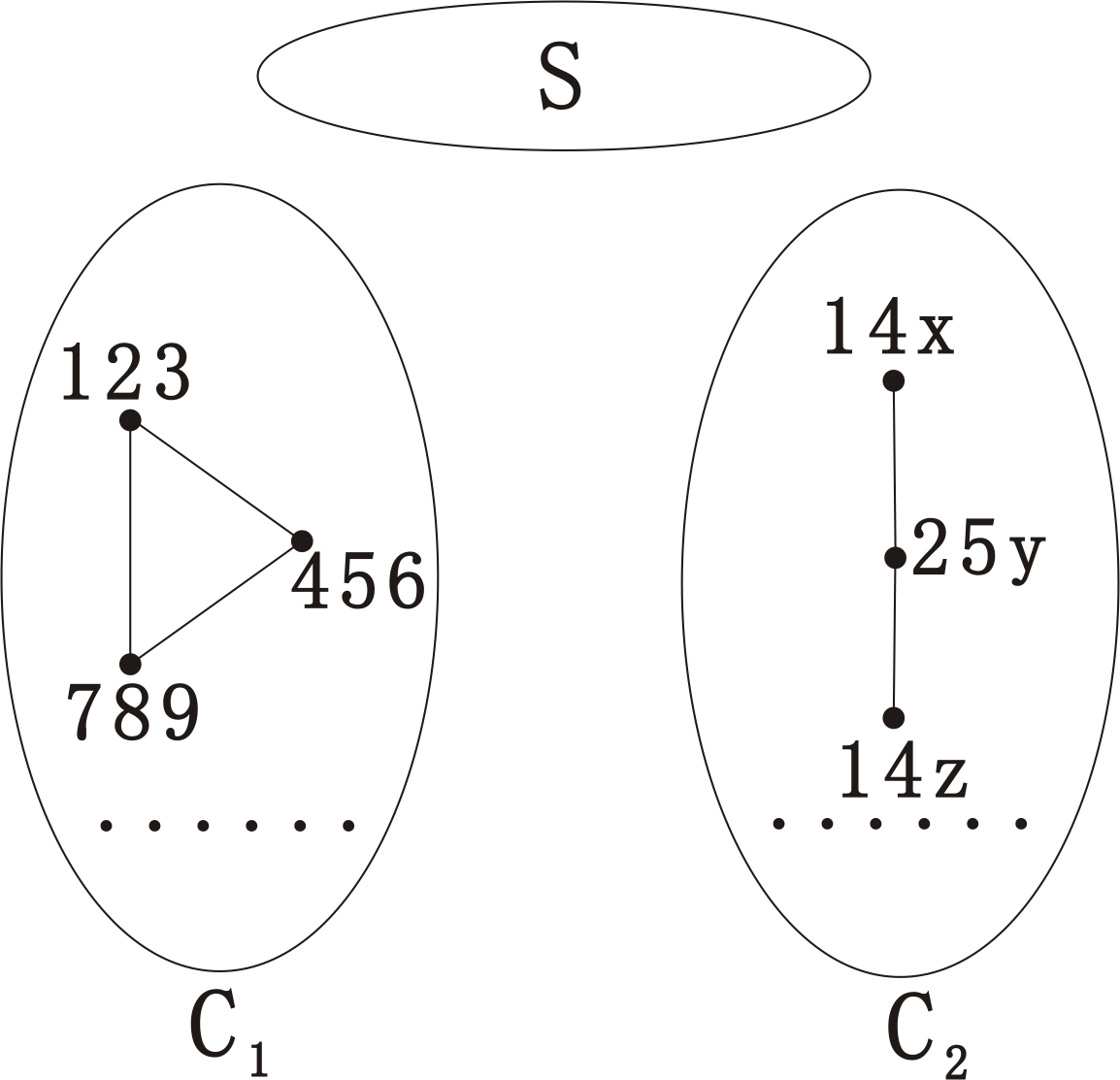}\\
  \caption{The case in $C_1$ and $C_2$}\label{case2.2}
\end{figure}

For the second case, assume the path consist of three vertices $\{1,4,x\}$, $\{2,5,y\}$, $\{1,4,z\} \in C_2$, where $x\neq y\neq z$ and $x,y,z \in \{7,8,9\}$. From the proof of Claim 3, we know that there are at most $6n-18$ vertices in $C_1$. Since we have double counted the $8$ vertices $\{1,4,y\}$, $\{1,5,7\}$, $\{1,5,8\}$, $\{1,5,9\}$, $\{1,6,y\}$, $\{2,4,7\}$, $\{2,4,8\}$, $\{2,4,9\}$. These vertices should be in either $C_1$ or $C_2$ but not in both. Thus, overall, $C_1 \cup C_2$ has no more than $18+6n-18-8=6n-8$ vertices, which is less than $9n-34$ when $n \geq 9$, then $C_1$ and $C_2$ has to be connected. A contradiction. See Figure \ref{case2.2} for an illustration.

Case 3: There is no $K_3$ in either $C_1$ or $C_2$, i.e. the components $C_1$ and $C_2$ contain $P_3$s. 

We shall consider the following three sub-cases based on the type of the paths.

Suppose there is a Type 1 path in $C_{1}$, let the three vertices be $v_{1}=\{1,2,3\}$, $v_{2}=\{4,5,6\}$, $v_{3}=\{1,7,8\}$, then by Claim 2, we know the number of vertices in $C_{2}$ is at most $3n+3$.

Now look at these vertices in $C_2$, there are at most $12$ vertices do not contain label $1$, if all of them are in $C_2$, i.e. none of them are included in the $S$, then these $12$ vertices form two cycles of length $6$. Of cause, if some of them are in $S$, then the rest of the vertices in each cycle form a set of paths. The rest of the vertices in $C_2$ all contain label $1$, thus not connected to each other, and they are connected to the vertices which do no contain label $1$. Then, we claim that either there is a Type 1 path, for example, $\{1,4,7\}$, $\{2,5,8\}$, $\{1,6,9\}$, or, there will be no more than $2n+4$ vertices in $C_2$. To see this, suppose we have no such desired path, and there are up to $n-2$ vertices containing both label $\{1,4\}$, there could be up to $n-2$ vertices in $C_2$ containing both labels $\{1,5\}$. Clearly not all $12$ vertices of no label $1$ are in $C_2$, since among those $12$ vertices, the ones such as $\{2,6,7\}$, $\{2,6,8\}$, $\{3,6,7\}$, $\{3,6,8\}$ will give us a desired path, thus there are at most $8$ among these $12$ vertices could be in $C_2$, also note that there must be some vertices from these $12$ vertices contained in $C_2$, otherwise, we have a set of singular vertices in $C_2$. Then the number of vertices in $C_2$ is at most $2(n-2)+8=2n+4$. If there are vertices containing both labels $\{1,6\}$ in $C_2$, then for sure we see the desired path. 

\begin{figure}[htbp]
  \centering
  \includegraphics[height=6cm]{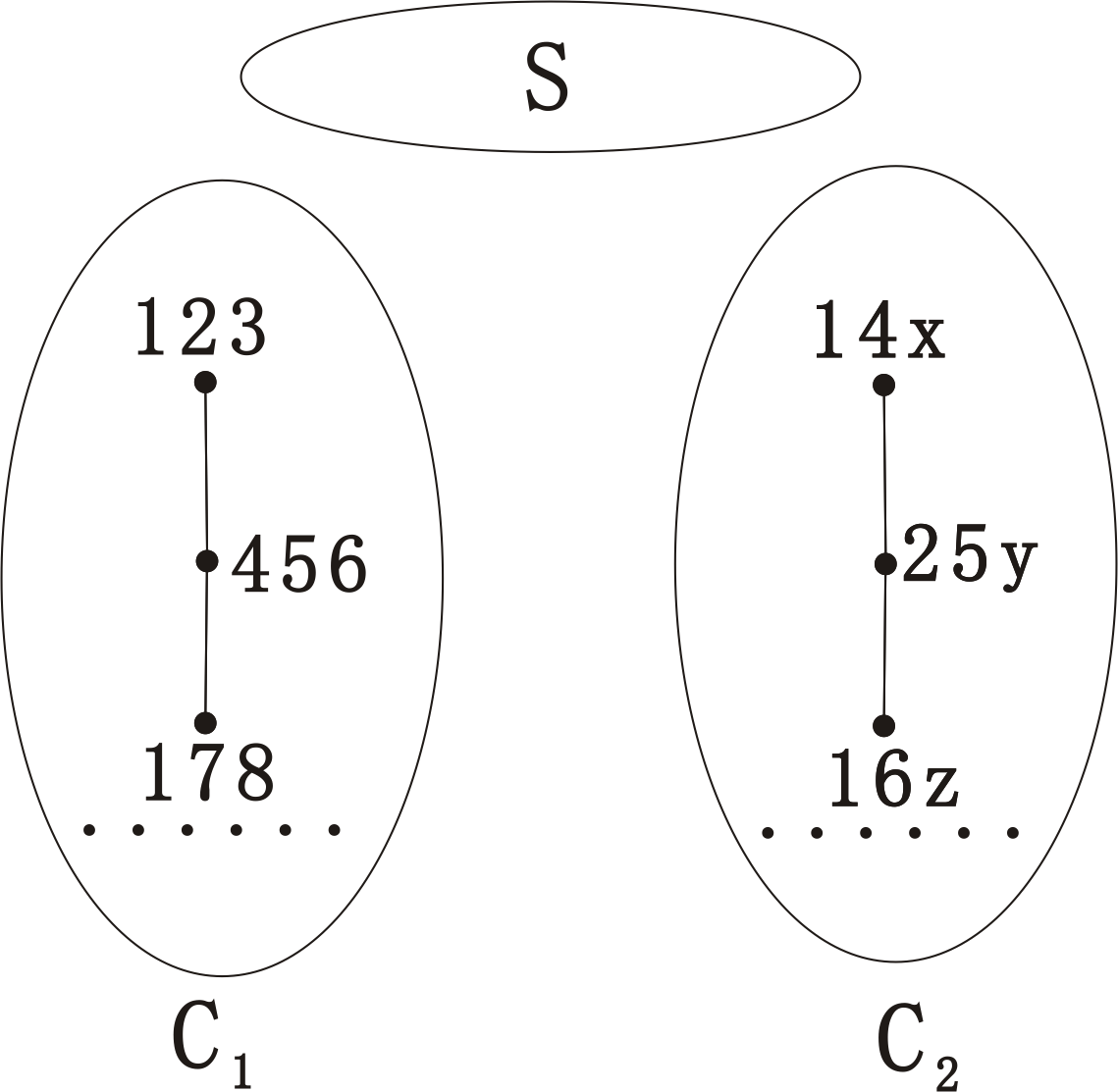}\\
  \caption{The case in $C_1$ and $C_2$}\label{case3.1}
\end{figure}

If we have the desired Type 1 path in $C_2$, let the three vertices be $\{1,4,x\}$, $\{2,5,y\}$, $\{1,6,z\}$, where $x\neq y\neq z$ and $x,z \in \{3,7,\cdots,n\}$, $y \in \{7,8\}$. Then based on the proof of Claim 2, $C_1$ has maximum $3n+3$ vertices, thus $C_1 \cup C_2$ has maximum $6n+6$ vertices. Also noticed that we have double counted the vertices of form $\{1,5,a\}$, where $a\in \{2,3,4,6,\cdots,n\}$, and vertices $\{1,2,4\}$, $\{1,2,6\}$, $\{1,4,y\}$, $\{1,6,y\}$, which both appear in the $C_1$ and $C_2$ in our calculation. Meanwhile, $\{2,4,6\}$ is only in $C_1$ or $S$, the vertex $\{3,5,7\}$ is either in $C_2$ or $S$. Depending on the choice of $x, y, z$, the vertex $\{3,5,7\}$ could also appear in $C_1$, for example, in the case $x=3,y=8,z=7$. If $\{3,5,7\}$ is either in $C_2$ or $S$, as $\{2,4,6\}$ and $\{3,5,7\}$ are connected, then one of them must be in $S$. If $\{3,5,7\}$ is in $C_1$, then $\{3,5,7\}$ is not in $C_1$, thus we know the size of $C_2$ has to be one less than the maximum possible.  The same for $\{1,2,7\}$ and $\{3,5,8\}$, $\{1,2,8\}$ and $\{3,6,7\}$. Therefore, there is no more than $5n+1$ vertices in $C_1 \cup C_2$, and $9n-34$ is larger than $5n+1$ when $n \geq 9$, then $C_1$ and $C_2$ has to be connected. A contradiction. See Figure \ref{case3.1} for an illustration.

\begin{figure}[htbp]
  \centering
  \includegraphics[height=6cm]{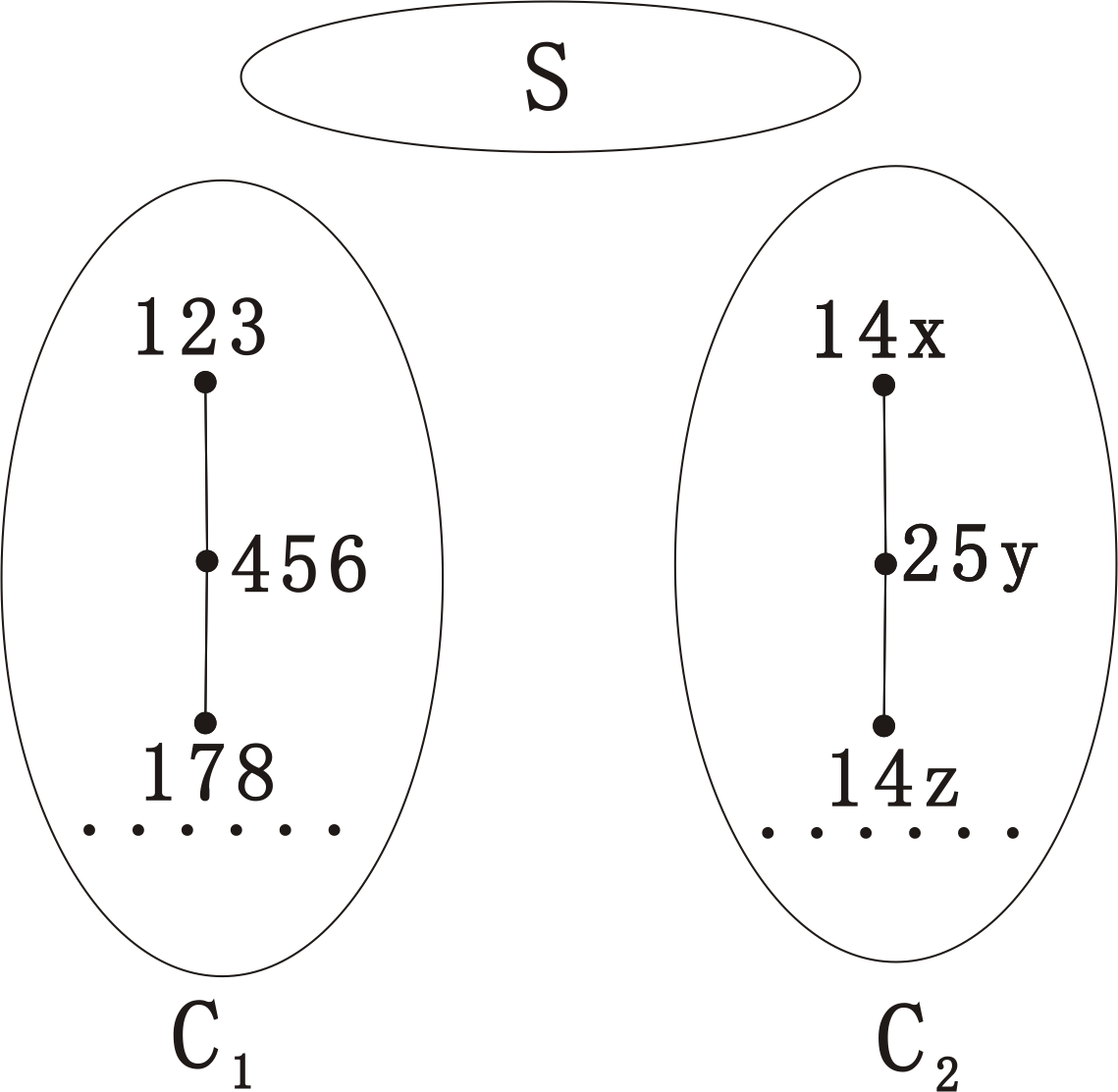}\\
  \caption{The case in $C_1$ and $C_2$}\label{case3.2}
\end{figure}

If there is no desired Type 1 path, then $C_2$ has at most $2n+4$ vertices, and we know there is a Type 2 path in $C_2$. Let the shared two labels be $\{1,4\}$ and let the three vertices be $\{1,4,x\}$, $\{2,5,y\}$, $\{1,4,z\}$ as show in Figure \ref{case3.2}, where $x\neq y\neq z$ and $x,z \in \{3,6,\cdots,n\}$, $y \in \{7,8\}$. Based on the proof of Claim 3, there are maximum $6n-18$ vertices in $C_1$. Since we have double counted the vertices of form $\{1,5,a\}$, where $a\in \{2,3,4,6,\cdots,n\}$, and vertices $\{1,2,4\}$, $\{1,2,6\}$, $\{1,4,y\}$, $\{1,6,y\}$, $\{2,4,7\}$, $\{2,4,8\}$, which both appear in the $C_1$ and $C_2$ in our calculation. Therefore, there is no more than $7n-18$ vertices in $C_1 \cup C_2$, and $9n-34$ is larger than $7n-18$ when $n \geq 9$, then $C_1$ and $C_2$ has to be connected. A contradiction.

Now assume that there is no Type 1 path in $C_1$, then there must be a Type 2 path in $C_1$. Let the three vertices on the path be  $v_{1}=\{1,2,3\}$, $v_{2}=\{4,5,6\}$, $v_{3}=\{1,2,7\}$, then by Claim 3, we know the number of vertices in $C_{2}$ is at most $6n-18$.


The case where there is a Type 2 path in $C_1$ and a Type 1 path in $C_2$ is the same to the case where there is a Type 1 path in $C_1$ and a Type 2 path in $C_2$. The latter we have considered above, so here we only consider the case where there is a Type 2 path in $C_1$ and there is also a Type 2 path in $C_2$.

\begin{figure}[htbp]
  \centering
  \includegraphics[height=6cm]{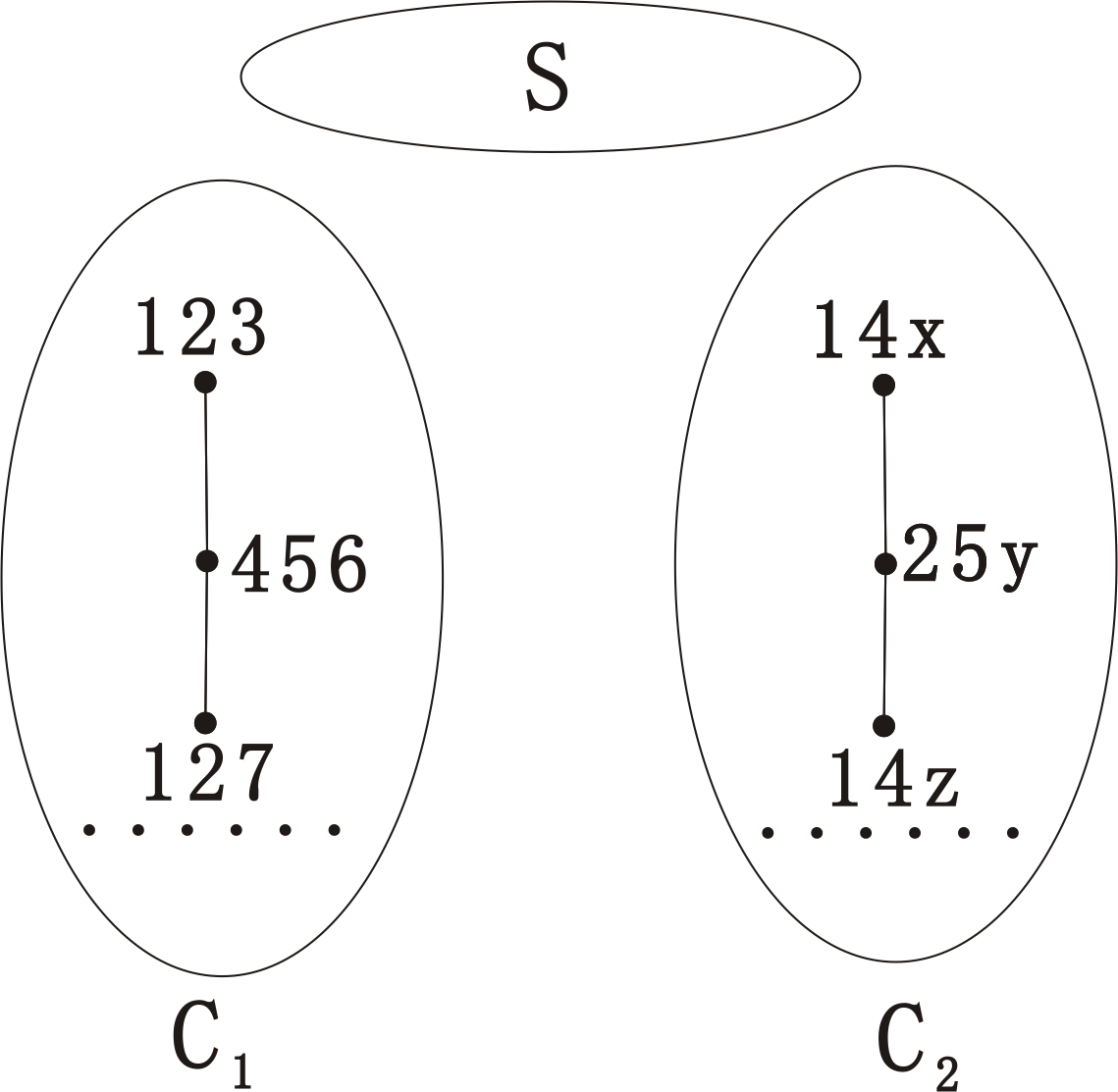}\\
  \caption{The case in $C_1$ and $C_2$}\label{case3.3}
\end{figure}

Suppose there are vertices containing both labels $\{1,4\}$ and vertices containing both labels $\{2,5\}$ in $C_2$, then there is no vertex containing both labels $\{1,6\}$ and no vertex containing both labels $\{2,6\}$ in $C_2$, and furthermore, there is no vertex containing neither label 1 nor label 2 in $C_2$, which implies that there is no vertex containing both labels $\{1,2\}$ in $C_2$, since the vertices with both $\{1,2\}$ only connect the vertices with no $\{1,2\}$ in $C_2$, otherwise Type 1 path or $K_3$ will appear in $C_2$. Then the number of vertices in $C_2$ is at most $4(n-3)-2$, i.e. at most $n-3$ vertices contain both labels $\{1,4\}$, at most $n-3$ vertices contain both labels $\{1,5\}$, at most $n-3$ vertices contain both labels $\{2,4\}$ and at most $n-3$ vertices contain both labels $\{2,5\}$, and we double counted the vertices $\{1,4,5\}$ and $\{2,4,5\}$. Then the Type 2 path in $C_2$ can be $\{1,4,x\}$, $\{2,5,y\}$, $\{1,4,z\}$, where $x\neq y\neq z$ and $x,y,z \in \{3,6,\cdots,n\}$, based on the prove of Claim 3, there have maximum $6n-18$ vertices in $C_1$. Since we have double counted the vertices of form $\{1,5,a\}$ and $\{2,4,b\}$, where $a\in \{2,3,4,6,\cdots,n\}$ and $b\in \{1,3,5,\cdots,n\}$, which both appear in the $C_1$ and $C_2$ in our calculation. Therefore, there is no more than $8n-28$ vertices in $C_1 \cup C_2$, and $G-S=9n-34$ is larger than $8n-28$ when $n \geq 9$, then $C_1$ and $C_2$ has to be connected. A contradiction. See Figure \ref{case3.3} for an illustration.

Thus we have showed that if $n\ge 9$, $G-S$ has to be connected, therefor the conjecture is true.

When $n=7$, only Type 2 path is possible in the graph. Let the three vertices in $C_1$ be $v_{1}=\{1,2,3\}$, $v_{2}=\{4,5,6\}$, $v_{3}=\{1,2,7\}$, then based on the proof of Claim 3, the number of vertices of $C_{2}$ is at most $6n-18=24$, that is, $9$ vertices contain label $1$, but not label $2$, $9$ vertices contain label $2$, but not label $1$, $3$ vertices contain both labels $\{1,2\}$ and $3$ vertices contain neither label $1$ nor label $2$.


As $C_2$ has at least three vertices, and it only has Type 2 path, let the three vertices on the path be $\{1,4,x\}$, $\{2,5,y\}$, $\{1,4,z\}$, where $x\neq y\neq z$ and $\{x,y,z\}=\{3,6,7\}$, then there are possibly three paths, it depends on the choice of $y$, the three paths are $\{1,4,3\}$, $\{2,5,6\}$, $\{1,4,7\}$ and $\{1,4,6\}$, $\{2,5,3\}$, $\{1,4,7\}$ and $\{1,4,3\}$, $\{2,5,7\}$, $\{1,4,6\}$. If the first path is presented in $C_2$, based on the proof of Claim 3, $C_1$ has at most $6n-18=24$ vertices, since we have double counted the $16$ vertices $\{1,2,4\}$, $\{1,2,5\}$, $\{1,2,6\}$, $\{1,3,5\}$, $\{1,3,6\}$, $\{1,4,5\}$, $\{1,4,6\}$, $\{1,5,6\}$, $\{1,5,7\}$, $\{1,6,7\}$, $\{2,3,4\}$, $\{2,4,5\}$, $\{2,4,6\}$, $\{2,4,7\}$, $\{3,5,7\}$, $\{3,6,7\}$. Meanwhile, $\{4,5,7\}$ can only be in $C_1$ or $S$, $\{2,3,6\}$ can only be in $C_2$ or $S$, but they are connected, then one of them must be in $S$, the same for $\{3,4,5\}$ and $\{2,6,7\}$, $\{3,4,6\}$ and $\{2,5,7\}$, $\{4,6,7\}$ and $\{2,3,5\}$. Thus, overall, no more than $24+24-16-4=28$ vertices, which is less than $|G-\kappa_1|=29$, then $C_1$ and $C_2$ has to be connected. A contradiction. We can discuss the second path and the third path in the same way, arrive the same conclusion, i.e. the number of vertices in $C_1 \cup C_2$ is at most $24+24-16-4=28$, which is less than $|G-\kappa_1|=29$, then $C_1$ and $C_2$ has to be connected. A contradiction.

When $n=8$, the three vertices in $C_{1}$ form a path $P_3$ of length $2$, it is possible for $C_1$ and $C_2$ to contain Type 1 path or Type 2 path, thus we have to look into each case. 

First, let $C_1$ has a Type 1 path $v_{1}=\{1,2,3\}$, $v_{2}=\{4,5,6\}$, $v_{3}=\{1,7,8\}$, then based on the proof of Claim 2, the number of vertices in $C_{2}$ is at most $3n+3=27$, that is, $15$ vertices contain label $1$ and $12$ vertices which do not label $1$.


If we have a Type 1 path in $C_2$, for example, $\{1,4,x\}$, $\{2,5,y\}$, $\{1,6,z\}$, where $x\neq y\neq z$ and $x,z \in \{3,7,8\}$, $y \in \{7,8\}$, then there are possibly four paths, they are $\{1,4,3\}$, $\{2,5,7\}$, $\{1,6,8\}$ and $\{1,4,8\}$, $\{2,5,7\}$, $\{1,6,3\}$ and $\{1,4,3\}$, $\{2,5,8\}$, $\{1,6,7\}$ and $\{1,4,7\}$, $\{2,5,8\}$, $\{1,6,3\}$, respectively. If the first path is presented in $C_2$, based on the proof of Claim 2, $C_1$ has at most $3n+3=27$ vertices, since we have double counted the $13$ vertices $\{1,2,4\}$, $\{1,2,5\}$, $\{1,2,6\}$, $\{1,3,5\}$, $\{1,4,5\}$, $\{1,4,7\}$, $\{1,5,6\}$, $\{1,5,7\}$, $\{1,5,8\}$, $\{1,6,7\}$, $\{2,4,8\}$, $\{3,5,8\}$, $\{3,6,7\}$. Meanwhile, $\{1,2,7\}$ can only be in $C_1$ or $S$, $\{3,4,8\}$ can only be in $C_2$ or $S$, but they are connected, then one of them must be in $S$, the same for the pairs $\{2,4,6\}$ and $\{3,5,7\}$, $\{4,5,8\}$ and $\{2,6,7\}$, $\{4,7,8\}$ and $\{1,3,6\}$. Thus, overall, no more than $27+27-13-4=37$ vertices, which is less than $|G-\kappa_1|=38$, then $C_1$ and $C_2$ has to be connected. A contradiction. We can discuss the second path, the third path and the fourth path in the same way and arrive the same contradiction.

If we have a Type 2 path in $C_2$, for example, $\{1,4,x\}$, $\{2,5,y\}$, $\{1,4,z\}$, where $x\neq y\neq z$ and $x,z \in \{3,6,7,8\}$, $y \in \{7,8\}$, then based on the proof of the Claim 3, $C_1$ has at most $6n-18=30$ vertices. Since we have double counted the $13$ vertices $\{1,2,4\}$, $\{1,2,5\}$, $\{1,2,6\}$, $\{1,3,5\}$, $\{1,4,5\}$, $\{1,4,y\}$, $\{1,5,6\}$, $\{1,5,7\}$, $\{1,5,8\}$, $\{1,6,y\}$, $\{2,4,7\}$, $\{2,4,8\}$, $\{3,4,y\}$. Meanwhile, $\{1,2,7\}$ can only be in $C_1$ or $S$. The vertex $\{3,5,8\}$ is either in $C_2$ or $S$. Depending on the choice of $x, y, z$, the vertex $\{3,5,8\}$ could also appear in $C_1$, for example, when $x=3,y=7,z=8$. If $\{3,5,8\}$ is either in $C_2$ or $S$, as $\{1,2,7\}$ and $\{3,5,8\}$ are connected, then one of them must be in $S$. If $\{3,5,8\}$ is in $C_1$, then we know the size of $C_2$ has to be one less than the maximum possible. 
The same for $\{1,2,8\}$ and $\{3,5,7\}$, $\{3,4,5\}$ and $\{2,6,7\}$, $\{4,5,7\}$ and $\{2,6,8\}$, $\{4,5,8\}$ and $\{3,6,7\}$, $\{2,4,5\}$ and $\{3,6,8\}$. $\{4,7,8\}$ can only be in $C_1$ or $S$, $\{1,3,6\}$ can only be in $C_2$ or $S$, but they are connected, then one of them must be in $S$. Thus, overall, no more than $27+30-13-7=37$ vertices, which is less than $|G-\kappa_1|=38$, then $C_1$ and $C_2$ has to be connected. A contradiction.

Second, let $C_1$ has a Type 2 path  $v_{1}=\{1,2,3\}$, $v_{2}=\{4,5,6\}$, $v_{3}=\{1,2,7\}$, then based on the proof of Claim 3, the number of vertices of $C_{2}$ is at most $6n-18=30$, that is, $9$ vertices contain label $1$, but not label $2$, $9$ vertices contain label $2$, but not label $1$, $3$ vertices contain both labels $\{1,2\}$ and $3$ vertices contain neither label $1$ nor label $2$.

The case where there is a Type 2 path in $C_1$ and a Type 1 path in $C_2$ is similar to the case where there is a Type 1 path in $C_1$ and a Type 2 path in $C_2$. The later we have considered already, so here we only consider the case where there is a Type 2 path in $C_1$ and there is also a Type 2 path in $C_2$.

Suppose, in $C_2$, there are vertices containing both labels $\{1,4\}$ and vertices containing both labels $\{2,5\}$ in $C_2$, then there is no vertex containing both labels $\{1,6\}$ in $C_2$, and there is no vertex containing both labels $\{2,6\}$ in $C_2$, otherwise Type 1 path will appear in $C_2$. Then the number of vertices in $C_2$ is at most $4\cdot5+6-2=24$, i.e. at most $5$ vertices contain both labels $\{1,4\}$, at most $5$ vertices contain both labels $\{1,5\}$, at most $5$ vertices contain both labels $\{2,4\}$ and at most $5$ vertices contain both labels $\{2,5\}$, at most $3$ vertices contain both labels $\{1,2\}$ and at most $3$ vertices contain neither label $1$ nor label $2$, and we double counted the vertices $\{1,4,5\}$ and $\{2,4,5\}$. Then the Type 2 path in $C_2$ can be $\{1,4,x\}$, $\{2,5,y\}$, $\{1,4,z\}$, where $x\neq y\neq z$ and $x,y,z \in \{3,6,7,8\}$, based on the prove of Claim 3, there have maximum $6n-18=30$ vertices in $C_1$. Since we have double counted the $14$ vertices $\{1,2,4\}$, $\{1,2,5\}$, $\{1,2,6\}$, $\{1,3,5\}$, $\{1,4,5\}$, $\{1,4,y\}$, $\{1,5,6\}$, $\{1,5,7\}$, $\{1,5,8\}$, $\{2,3,4\}$, $\{2,4,5\}$, $\{2,4,6\}$, $\{2,4,7\}$, $\{2,4,8\}$, which both appear in the $C_1$ and $C_2$ in our calculation. Meanwhile, $\{3,4,5\}$ can only be in $C_1$ or $S$, the vertex $\{1,6,7\}$ is either in $C_2$ or $S$. Depending on the choice of $x, y, z$, the vertex $\{1,6,7\}$ could also appear in $C_1$, for example, when $x=6,y=7,z=8$. If $\{1,6,7\}$ is either in $C_2$ or $S$, as $\{3,4,5\}$ and $\{1,6,7\}$ are connected, then one of them must be in $S$. If $\{1,6,7\}$ is in $C_1$, then we know the size of $C_2$ has to be one less than the maximum possible. The same for pairs $\{4,5,7\}$ and $\{1,6,8\}$, $\{4,5,8\}$ and $\{1,3,6\}$, $\{1,2,8\}$ and $\{3,5,7\}$. Therefore, there is no more than $24+30-14-4=36$ vertices in $C_1 \cup C_2$, and $|G-\kappa_1|=38$ is larger than $36$, then $C_1$ and $C_2$ has to be connected. A contradiction.

In summary, we have proved that when $k=3$, the conjecture is true and the bound is achieved only in the case that one of the disconnected component contains just two vertices linked by an edge.

\end{document}